\documentclass[11pt]{article}
\usepackage{amsmath,amsthm,amssymb}
\usepackage[applemac]{inputenc}
\usepackage{mathrsfs}
\usepackage{graphicx}
\usepackage{amsfonts}
\theoremstyle{plain}

\theoremstyle{definition}

\numberwithin{equation}{section}

\newcommand{\pical}{\mathcal{P}}
\newcommand{\tto}{\to}

\newcommand{\haus}{\mathcal H}
\newcommand{\huno}{\haus^1}

\newcommand{\M}{\mathcal M}

\newcommand{\lcal}{\mathcal L}

\newcommand{\F}{\mathcal F}

\newcommand{\cH}{\mathcal H}
\newcommand{\W}{\mathcal W}

\newcommand{\si}{\sigma}

\newcommand{\Om}{\Omega}

\newcommand{\ve}{\varepsilon}

\newcommand{\N}{\mathbb N}

\newcommand{\R}{\mathbb R}

\title{Models and applications of Optimal Transport\\ in Economics, Traffic and Urban Planning}
\author{Filippo Santambrogio\thanks{\scriptsize\ CEREMADE, UMR CNRS 7534, Universit\'e Paris-Dauphine, Pl. de Lattre de Tassigny, 75775 Paris Cedex 16, FRANCE
\texttt{filippo@ceremade.dauphine.fr , http://www.ceremade.dauphine.fr/$\sim$filippo}.}}
\date{Grenoble, June 2009\\ Revised version : March 23rd, 2010}
\begin{document}
\maketitle

\begin{abstract}
Some optimization or equilibrium problems involving somehow the concept of optimal transport are presented in these notes, mainly devoted to applications to economic and game theory settings. A variant model of transport, taking into account traffic congestion effects is the first topic, and it shows various links with Monge-Kantorovich theory and PDEs. Then, two models for urban planning are introduced. The last section is devoted to two problems from economics and their translation in the language of optimal transport.
\end{abstract}

\medskip\noindent
{\bf AMS Subject Classification (2010):} 00-02, 90B06, 49J45, 90C25, 91B40, 91B50, 91D10 35Q91 

\bigskip\noindent
{\bf Keywords:} traffic congestion, wardrop equilibrium, optimal flows, equilibrium problems, Kantorovich potential, contract theory, Nash equilibrium

\tableofcontents

\section*{Introduction}

These lecture notes will present the main issues and ideas of some variational problems that use or touch the theory of Optimal Transportation. They all come from economic-oriented applications. Problems will be presented through the main ideas, with almost no proofs. As I did during the class in Grenoble and for the other lecture notes (\cite{GrenobleIntro}), I will try to keep a very informal level of presentation.

The first topic I will present is a variant of the usual Monge problem, which takes into account congestion effects in the transportation. It has a more dynamical taste since it asks for looking at the trajectories followed by the particles instead of simple pairs $(x,T(x))$. The starting point will be the equivalent formulation by Beckmann, modified in order to take into account non-uniform metrics and then traffic intensity. In this way one obtains a model which is linked to game theory and Nash equilibria and has a well-known discrete counterpart on networks. Moreover, the optimality conditions for the minimization let 
a classical transport problem appear (for a metric which is not known a priori).

After this first section, two models for the distribution of some fundamental elements of the structure of urban regions (residents, jobs, industrial areas, services\dots) are discussed. The first is purely variational: we suppose that these distribution optimize a total welfare functional, as if a benevolent planner could control the city; the functional let a transport cost (say, a Wasserstein distance) appear explicitly. The second deals with much more delicate equilibrium issues on the agents' choices and it does not addresses explicitly any transport problem in its formulation; yet, the theory of Monge-Kantorovitch is very useful in its resolution. In both cases the Kantorovitch potential play an essential role.

By the end of this third section, the reader who is not familiar with economic theories should have got more accustomed to the some typical concepts of the rational behavior of the consumers, and he will be ready for Section 4. This section presents two classical problems in economics (in situations of competition and of monopoly, respectively) and it shows how to translate them into problems involving transport costs and Kantorovitch potentials.

 \section{Traffic congestion}
\subsection{Generalizations of Beckmann's Problem}

We saw in the introductory lecture notes \cite{GrenobleIntro} the problem (B):
$$
 \min \left\{M(\lambda)\,: \lambda\in\M^d(\Omega); \;\nabla\cdot\lambda=\mu-\nu\right\},
$$
where $M(\lambda)$ denotes the mass of the vector measure $\lambda$ and we said that it is equivalent to the original problem of Monge. Actually, one way to produce a solution to this divergence-constrained problem, is the following:  take an optimal transport plan $\gamma$ and build a vector measure $v_\gamma$ defined through
$$
<v_\gamma,\phi>:=\int_{\Omega\times\Omega}\int_0^1 \omega'_{x,y}(t)\cdot \phi(\omega_{x,y}(t))dt\, d\gamma,
$$
for every $\phi\in C^0(\Omega;\R^d)$, $\omega_{x,y}$ being a parametrization of the segment $[x,y]$.

It is not difficult to check that this measure satisfies the divergence constraint, since if one takes $\phi=\nabla \psi$ then 
$$\int_0^1 \omega'_{x,y}(t)\cdot \phi(\omega_{x,y}(t))=\int_0^1 \frac{d}{dt} \left(\psi(\omega_{x,y}(t)\right)dt=\psi(y)-\psi(x)$$
and hence  $<v_\gamma,\nabla\psi>=\int\psi\,d(\nu-\mu)$.

To estimate its mass we can see that $|v_\gamma|\leq \sigma_\gamma$, where the scalar measure $\sigma_\gamma$ is defined through
$$
<\sigma_\gamma,\phi> :=\int_{\Omega\times\Omega}\int_0^1 |\omega'_{x,y}(t)|\phi(\omega_{x,y}(t))dt\, d\gamma,\quad\forall\phi\in C^0(\Omega;\R)
$$
and it is called {\it transport density}. The mass of $\sigma_\gamma$ is obviously
$$\int d\sigma_\gamma=\int\int_0^1 |\omega'_{x,y}(t)|dt\, d\gamma=\int |x-y|d\gamma=W_1(\mu,\nu),
$$
which proves the optimality of $v_\gamma$.

It is interesting to investigate whether  $\sigma_\gamma<<\lcal^d$, since this would imply that Problem (B) is well-posed in $L^1$ instead of the space of vector measure. For the sake of the variants that we will see later on, it would be interesting to give conditions so that $\sigma_\gamma\in L^p$ as well. All these subjects have been widely studied by De Pascale, Pratelli (see \cite{DePPra, DePEvaPra, DePPra2}) but there is a more recent (and shorter) proof of the same estimates in \cite{simple proof}. It is in particular true that $\mu,\nu\in L^p$ implies that $\sigma_\gamma\in L^p$ and that it is sufficient that one of the two measures is absolutely continuous in order to get the same on $\sigma_\gamma$.

The simplest possible generalization of Problem (B) is the following:  
  $$\min \int k(x)|v(x)|dx\;:\;\nabla\cdot v =\mu-\nu$$
  that corresponds, by duality with the functions $u$ such that $|\nabla u|\leq k$, to 
     $$\min \int d_k(x,y)d\gamma\;:\;\gamma\in\Pi(\mu,\nu),$$
 where $d_k(x,y)=\inf_{\omega(0)=x,\,\omega(1)=y}L_k(\omega):=\int_0^1 k(\omega(t))|\omega'(t)|dt$ is the distance associated to the Riemannian metric $k$. It would be possible to build in this case an optimal $v_\gamma$ by replacing the curves $\omega_{x,y}$ with the $k-$geodesics (instead of the segments). 

This generalization above comes from the modelization of a non-uniform cost for the movement (due to geographical obstacles or configurations). It can be applied to several situation but it is anyway evident that one should look for more realistic models, at least in the case of urban transport. In this case the metric $k$ is usually not a priori known, but it depends on the traffic distribution itself.

The simplest model could be considering a metric $k(x)=g(|v(x)|)$ depending through an increasing function $g$ on the traffic itself (represented by the intensity of $v$). In this case a very naive model would be obtained by setting $H(t)=tg(t)$ and then solving
$$\min \int H(|v(x)|)dx\;:\;\nabla\cdot v =\mu-\nu.$$
In most cases, $H$ is strictly convex and this is a strictly convex counterpart to the problem by Beckmann (which was somehow suggested by Beckmann himself in his book \cite{BecPuu}). Notice that this model is not completely realistic neither since it allows for ``cancellation'' effects: several flows in opposite directions at a same point $x$ may give a total vector $v(x)=0$, even if the number of travellers at $x$ is high. Yet, in Section 3.3 we will see that this simplifed model will turn out to be equivalent to a more precise one. 

We just mention that there exist concave variant too, which are known under the name of {\it branched transport}.  This name is used for addressing all the transport problems where the cost for a mass $m$ moving on a distance $l$ is proportional to $l$ but subadditive w.r.t. $m$. Typically, it is proportional to a power $m^\alpha$ ($0<\alpha<1$).  The adjective ``branched'' in the name stands for one of the main features of the optimal solutions: they gather mass together, masses tend to move jointly as long as possible, and then they branch towards different destinations, thus giving rise to a tree-shaped structure.

This problem comes from a discrete problem on graphs, where the cost of a graph $G$ whose edges  $e_h$ are weighted with coefficients $w_h$ is of the form $\sum_h p_h^\alpha \huno(e_h)$. It has a continuous generalization where the energy to be minimized is
$$M^\alpha(v)=\begin{cases} \int_M \theta^\alpha d\huno &\mbox{ if } v = U(M,\theta,\xi),\\
                                                                +\infty &\mbox{ otherwise. } \end{cases}
$$
     where $v = U(M,\theta,\xi)$ means that $v$ is a rectifiable measure supported on the set $M$, with orientation $\xi$ and density (multiplicity) $\theta$. The energy $M^\alpha$ is then minimized under the constraint $\nabla\cdot v =\mu-\nu$.
     
     These notes will not develop any more this alternative problem and the reader may find the whole theory of branched transport in the recent book by Bernot, Morel and Caselles (\cite{Book Irrigation}).
     
     \subsection{Wardrop equilibria, the discrete case}
     We will describe in this section a traffic problem which has some interesting issues on equilibria and some interesting relations with optimal transport theory. We will start from the discrete case on networks and then generalize to the continuous case. The network case was introduced in \cite{wardrop} and then studied in \cite{BecMcGWin}.
     
In the discrete framework, one considers

  \begin{itemize}
		\item A finite graph with edges $e\in E$ and a set of sources $S$ and destinations $D$, 
		\item the set $C(s,d)=\{\omega\mbox{ from } s \mbox{ to }d\}$ of possible paths from $s$ to $d$, 
		\item a demand input $\gamma=(\gamma(s,d))_{s,d}$ denoting the quantity of commuters from each $s\in S$ to each $d\in D$, or a set $\Gamma$ of possible $\gamma$ (for instance this could be the set of all demands where the total number of commuters leaving each point $s$ and the total number arriving to each point $d$ are prescribed, but the coupling, i.e. how many commuters for each pair $(s,d)$, is not);  
		\item an unknown repartition strategy (to be looked for) $q=(q_\omega)_\omega$ such that $\sum_{\omega\in C(s,d)}q_\omega=\gamma(s,d)$, 
		\item a consequent traffic intensity on each edge $e$ (depending on $q$) $i_q=(i_q(e))_e$ given by $i_q(e)=\sum_{e\in\omega}q_\omega$, 
		\item an increasing function $g:\R^+\to\R^+$ such that $g(i_q(e))$ represents the congestioned cost of $e$, 
		\item the cost for each path $\omega$, given by $c(\omega)=\sum_{e\in\omega}g(i_q(e))$.
	\end{itemize} 

  The global strategy $q$ represents the overall distribution on choices of commuters' paths. Imposing a {\it Nash equilibrium} condition (no single commuter wants to change his choice, provided all the others keep the same strategy) gives the following condition:
  $$\omega\in C(s,d),\,q_\omega>0\Longrightarrow c(\omega)=\min\{c(\tilde{\omega})\,:\,\tilde{\omega}\in C(s,d)\}.$$ 
  This condition is well-known among geographical economists as {\it Wardrop equilibrium}. 
  
 The existence of at least an equilibrium comes from the following variational principle.
 
Optimizing an overall congestion cost means minimizing a quantity $J(q):=\sum_e H(i_q(e))$ (where $H:\R^+\to\R^+$ is an increasing function: for instance if one takes $H(t)=tg(t)$ the value of $J(q)$ gives the total cost for all commuters) among all possible strategies $q$.   

The minimization of $J$ has obviously a solution and one can look for  optimality conditions. Suppose that $H$ and $\Gamma$ are convex, so that the necessary conditions will also be sufficient: it is easy to see that $q$ minimizes if and only if, for every other admissible $\tilde{q}$, one has
$$\sum_e H'(i_q(e))(i_{\tilde{q}}(e)-i_q(e))\geq 0.$$
Set $\xi(e):=H'(i_q(e)$ and rewrite the right hand side has
$$\sum_e \xi(e)(i_{\tilde{q}}(e)-i_q(e))=\sum_e\sum_{\omega\ni e}\xi(e)(\tilde{q}(\omega)-q(\omega))=\sum_\omega\left(\sum_{e\in \omega} \xi(e)\right)(\tilde{q}(\omega)-q(\omega)).$$
This says that, if one sets $L_\xi(\omega):=\sum_{e\in \omega} \xi(e),$ the optimal $q$ must minimize $\sum_\omega L_\xi(\omega)q(\omega),$ since we got $\sum_\omega L_\xi(\omega)\tilde{q}(\omega)\geq  \sum_\omega L_\xi(\omega)q(\omega)$.

This means two facts. First, since the conditions of admissibility on $q$ only look at starting and arrival points, it is pointless to put some mass on those curves $\omega$ from a source $s$ to a destination $d$ such that the value $L_\xi(\omega)$ is strictly larger than $d_\xi(s,d):=\min_{\omega\in C(s,d)}L_\xi(\omega)$. This means that $q(\omega)>0$ and $\omega\in C(s,d)$ imply $L_\xi(\omega)=d_\xi(s,d)$.

Second, another condition occurs when the demand $\gamma$ is not fixed.  We said that to optimize  $\sum_\omega L_\xi(\omega)q(\omega)$ we only use curves where $L_\xi=d_\xi$ and this gives 
$$\sum_\omega L_\xi(\omega)q(\omega)=\sum_{s,d}d_\xi(s,d)\left(\sum_{\omega\in C(s,d)}q_\omega\right)=\sum_{s,d}d_\xi(s,d)\gamma(s,d).$$
In particular, one also needs to choose $\gamma\in\Gamma$ so as to minimize
$$\sum_{s,d}d_\xi(s,d)\gamma(s,d),\quad \gamma\in\Gamma.$$
This second condition is empty if $\Gamma$ only contains one $\gamma$ but it is of particular interest when $\Gamma=\Pi(\mu,\nu)$, since it says that $\gamma$ must solve a Kantorovitch problem for the cost $d_\xi$.

The first condition, on the other hand, always gives some information on $q$ and exactly says:  if $q$ is optimal, then it is a Wardrop equilibrium for $g=H'$.

     \subsection{Wardrop equilibria, the continuous case and equivalences}

It is possible to give a continuous formulation and prove analogous results (see \cite{CarJimSan}).
  In a domain $\Omega\subset\R^n$ the demands are represented by probabilities $\gamma\in\pical(\Omega\times\Omega)$. We are given a set $\Gamma\subset\pical(\Omega\times\Omega)$ as the set of admissible demand couplings:   usually $\Gamma=\{\bar{\gamma}\}$ or $$\Gamma=\Pi(\mu,\nu)=\{\gamma\in\pical(\Omega\times\Omega)\,:\,(\pi_X)_\sharp\gamma=\mu,\,(\pi_Y)_\sharp\gamma=\nu\}.$$ 
  Let us also set
\begin{gather*}
C=\{\mbox{Lipschitz paths } \omega:[0,1]\to\Omega\}\\
C(s,d)=\{\omega\in C: \omega(0)=s,_,\omega(1)=d\}.
\end{gather*} 
  We look for a probability $Q\in\pical(C)$ such that $(\pi_{0,1})_\sharp Q\in\Gamma$.
  
We want to define a traffic intensity $i_Q\in\mathcal{M}^+(\Omega)$ such that the quantity $i_Q(A)$ stands for ``how much '' the movement takes place in $A$\dots  For $\phi\in C^0(\Omega)$ and $\omega\in C$ set  $L_\phi(\omega)=\int_0^1\phi(\omega(t)|\omega'(t)|dt.$
  
  Then we define $i_Q$ by
  $$<i_Q,\phi>=\int_C L_\phi(\omega)Q(d\omega)=\int_C\Big( \int_0^1 \varphi(\omega(t)) \vert \omega'(t))\vert dt \Big)Q(d\omega).$$ 
  
  Notice that this is exactly what happens for the transport density! the traffic intensity $i_Q$ is a generalization of the transport density, since it deals with the case where $Q$ is any measure on $C$, while the transport density only looks at the measure concentrated on the segments $[x,y]$
for $(x,y)$ in the support of an optimal $\gamma$. 

 In this continuous framework, it is more convenient to start from the optimization point of view (instead of looking at the equilibrium as a starting point): we minimize the convex functional
  $$J(Q)=\begin{cases}\int H(i_Q(x))dx&\mbox{ if }i_Q<<\lcal^n,\\
                        +\infty         &\mbox{otherwise}\end{cases}$$
 among all admissible strategies $Q$, $H$ being a convex, increasing and superlinear function. Typically $H(t)=t^p$, or $H(t)=t+t^p$ (which is more reasonable since in general we have $g(0)=H'(0)$ and we do not want $g(0)=0$: this would mean that moving on an empty road costs nothing, which is usually not the case).
 
 First one should prove finiteness of the minimum, which is not evident since in the continuus case one needs to prove the existence of a $Q$ such that $i_Q\in L^p$. This is, in the case of $\mu,\nu\in L^p$, a consequence of the summability results on the transport density, since the transport density is, as we said, a particular choice for $i_Q$. This is why we explicitly cited the $L^p$ result fo De Pascale and Pratelli (besides its interest in itself).

It is possible to look for optimality conditions and to reobtain the same Wardrop equilibrium + optimization of $d_\xi$. Here $\xi$ will be the metric $\xi(x)=H'(i_Q(x))$. Yet, this function is not continuous nor l.s.c. and some efforts should be spent to give a meaning to the concept of geodesic distance in the case $\xi\in L^q$.

It is also interesting to notice that this problem looks at the movement of some players whose individual goal is fixed but whose utility also looks at the density of all the other players (i.e. their movement is more expensive if they pass where the density is higher): this seems to be a particular case of the so-called {\it Mean Field Games} introduced by Lasry and Lions in \cite{LasLio}.

All the results we cited are valid for the case $\Gamma=\{\bar{\gamma}\}$ as well as for $\Gamma=\Pi(\mu,\nu)$ (all the transport plans).

Yet, in this second case, something more may be said. 
Instead of defining a scalar traffic intensity one can define a vector measure $v_Q$ by:
   \[\int_{\overline{\Omega}} \varphi(x)\, d v_Q(x):=\int_{C([0,1];\overline{\Om})}\Big( \int_0^1 \varphi(\omega(t)) \cdot  \omega'(t)) dt \Big)dQ(\omega) , \; \forall \varphi \in C(\overline{\Omega}, \R^N),\]
   i.e. sort of a vector version of $i_Q$.
It is immediate to check that $\vert v_Q \vert \leq i_Q$, and that
\[\nabla\cdot v_Q=\mu-\nu, \; \ v_Q\cdot\nu=0 \mbox{ on } \partial\Om.\]
Since $H$ is increasing, this implies that the infimum of the previous problem with $i_Q$ is larger than that of the minimal flow problem: 
\begin{equation}
\label{miniflow}
\inf\left\{\int_\Om \cH(v)\ dx\ :\ \nabla\cdot v=\mu-\nu,\ v\cdot\nu=0 \mbox{ on } \partial\Om\right\},
\end{equation}
where $\cH(v):=H(\vert v\vert)$. 

A natural question, arising for instance from a comparison with the Monge case, where looking for the vector or the scalar transport density was the same, is the possible equivalence of the two problems. 

One can see that a minimizer of the scalar problem can be built formally from a minimizer of the vector one in the following way:
if $v$ is the unique solution of the vector problem \eqref{miniflow} and $\mu$ and $\nu$ are absolutely continuous (so that we will write $\mu$ and $\nu$ for their densities as well), we consider the non-autonomous Cauchy problem 
\begin{equation}
\label{problemacauchy}
\left\{\begin{array}{ccc}
\omega'(s)&=&w(s,\omega(s))\\
\omega(0)&=&x\\
\end{array}            \right.       {}
\end{equation}
for the non-autonomous vector field
\begin{equation}
\label{modified}
w(t,x)=\frac{v(x)}{(1-t)f_0(x)+t f_1(x)},\quad (t,x)\in[0,1]\times\Om.
\end{equation}
The latter will not have
any Lipschitz continuity property in general, unless the optimizer $v$ of \eqref{miniflow} is regular: anyway, if we assume that one can prove $v\in \mathrm{Lip}(\Om)$, then the flow $X:[0,1]\times\Om\to\Om$ of $v$ is well-defined as the solution of \eqref{problemacauchy} and we can take $\mu_t$ as the image of $\mu$ through the map $X(t,\cdot)$. One can see that $\mu_t$ must coincide with  the linear interpolating curve
$(1-t)\mu+t\nu$ (because this curve solves the continuity equation thanks to the divergence condition).
This yields that
$
(X(1,\cdot))_\sharp f_0=f_1,
$
which ensures that $X(1,\cdot)$ transports $\mu$ on $\nu$. If we now consider the probability measure concentrated
on the flow, i.e. 
\[
Q=\delta_{X(\cdot,x)}\otimes \mu,
\]
then $Q$ is admissible and it is not difficult to see that $i_Q=|\si|$, since 
$$\int \phi d i_Q=\int \int_0^1 \phi(\omega_x(t))|\omega_x'(t)|\,dt \,d\mu=\int_0^1 dt \int \phi(\omega_x(t))\frac{|v(\omega_x(t))|}{\mu_t(\omega_x(t))} \,d\mu=\int_0^1 dt \int \phi\frac{|v|}{\mu_t} \,d\mu_t=\int \phi |v|.$$ 
This finally implies that the minima of the two problems coincide. Moreover, this construction provides a transport map (that is $X(1,\cdot)$) from $\mu$ to $\nu$, whose transport
``rays'' evidently do not cross and which is monotone on transport ``rays'' (as a consequence of Cauchy-Lipschitz Theorem). 

Notice that if one wanted to prove rigorously what we stated he should investigate a little bit the regularity of the optimal $v$. This may be done if one writes optimality conditions for $v$ and sees that he has $v=\nabla \cH^*(\nabla u)$ where $u$ solves 
\begin{equation}
\label{euler}
\left\{\begin{array}{cccc}
\nabla\cdot\nabla\cH^*(\nabla u)&=&\mu-\nu,& \mbox{ in } \Om,\\
\nabla\cH^*(\nabla u)\cdot \hat{n}&=&0, & \mbox{ on } \partial\Om,
\end{array}
\right.
\end{equation}

For $H(t)=t^2$ this is a simple Laplace equation and regularity theory is well-known. For $H(t)=t^p$ this gives a $p'-$Laplace equation and here as well lots of studies have benn done. Yer, for modeling reasons, we said that it is important to look at the case $H'(0)>0$, and we suggested as a typical case \begin{equation}
\label{H}
\cH(\si)=\frac{1}{p}|v|^p+a|v|,\ v\in\R^N,
\end{equation}
which leads to a function $\cH^*$ which vanishes on $\overline{B_1}$. In particular, the corresponding equation for $u$ is very very degenerate and regularity results are less studied (see \cite{BraCarSan}, both for the equivalence with the Wardrop problem and for some regularity proofs).

 \section{The urban planning of residents and services}
 
 A very simplified model that has been proposed for studying the distribution of residents and services in a given urban region $\Omega$ passes through the minimization of a total quantity $\F(\mu,\nu)$ concerning two unknown densities $\mu$ and $\nu$.

The two measures $\mu$ and $\nu$ will be searched among probabilities on $\Omega$. This 
means that the total amounts of population and production are 
fixed as problem data. The definition of
the total cost functional to optimize takes into account some criteria we want the two
densities $\mu$ and $\nu$ to satisfy:

\begin{itemize}
\item[(i)] there is a transportation cost $T$ for moving from the residential
areas to the services areas;

\item[(ii)] people do not want to live in areas where the density of
population is too high;

\item[(iii)] services need to be concentrated as much as possible in order
to increase efficiency and decrease management costs.

\end{itemize}

Fact (i) is described, in its easiest version, through a $p$-Wasserstein
distance ($p\geq1$). We will look at $T(\mu,\nu)=W_p^p(\mu,\nu)$.

Fact (ii) will be described by a penalization functional, a kind of total
unhappiness of citizens due to high density of population, obtained by integrating
with respect to the citizens' density their personal unhappiness. 

Fact (iii) is modeled by a third term representing costs for managing services once
they are located according to the distribution $\nu$, taking into account that
efficiency depends strongly on how much $\nu$ is concentrated.

The cost functional to be considered is then
\begin{equation}\label{e1costf}
\F(\mu,\nu)=T(\mu,\nu)+F(\mu)+G(\nu),
\end{equation}
where $F,G:\pical(\Omega)\to [0,+\infty]$ are functionals chosen so that the first one favors spread measures and the second one concentrated measures, in suitable senses.

We stress that this model is a very naive one, since it disregards equilibrium issues and several other parameters, and that it could be applied only in those cases where a planner could control the whole behavior of the region. We refer to \cite{3opti, ButSan, ButSanSIREV, T+GV, PhDThesis} for the study of this model and of similar ones.

As far as particular choices for the functionals $F$ and $G$ are concerned, we may consider
\begin{gather*}
F(\mu)=\begin{cases}\int_{\Omega}f(u)\,d\lcal^d&\text{ if }\mu=u\cdot\lcal^d\\ +\infty&\text{ otherwise,}\end{cases}\\
G(\nu) =
\begin{cases} \sum_{k\in\N}g(a_k) & \textrm{ if } \nu=\sum_{k\in\N}a_k\delta_{x_k} \\ +\infty & \textrm{otherwise,}\end{cases}
\end{gather*}
where the integrand $f:[0,+\infty]\tto[0,+\infty]$ is assumed to
be lower semicontinuous and convex, with $f(0)=0$ and superlinear at
infinity, that is,
$$\label{superlinearita}\lim_{t\tto+\infty}\frac{f(t)}{t}=+\infty,$$
and the function $g$ is required to be subadditive, lower semicontinuous, and such that 
$$g(0)=0\quad\mbox{ and }\lim_{t\tto 0}\frac{g(t)}{t}=+\infty.$$

In this form we have two {\it local} lower semicontinuous functional on measures (see \cite{bb1}: a functional on measures is said to be local if it is additive on mutually singular measures ). This is a useful class of functionals over measures including both
concentration preferring functionals and functionals
favoring spread measures. 

Without loss of generality, 
by subtracting constants to the functional $F$, we can suppose $f'(0)=0$.
Due to the assumption $f(0)=0$, the ratio $f(t)/t$ is an incremental ratio 
of the convex function $f$ and thus it is increasing in $t$. Then, if we write the functional $F$ as
$$\int_{\Omega}\frac{f(u(x))}{u(x)}u(x)\,dx,$$
we can see the quantity $f(u)/u$, which is increasing in $u$, as the unhappiness of a single 
citizen when he lives in a place where the population density is $u$. Integrating it with respect 
to $\mu=u\cdot\lcal^n$ gives a quantity to be seen as the total unhappiness of the population. 

Concerning $G$, we can think that we are requiring $\nu$ to be concentrated on a limited number of service poles and that the effects of the managing costs and of the production of a pole whose size is $a$ are summarized in a cost function $g(a)$. 

For $G$, there are other interesting choices among functionals which favor concentration. One of them could be
$$G(\nu)=\int_\Omega \int_\Omega h(|x-y|)\,\nu(dx)\nu(dy),$$
where $h$ is an increasing function and $h(|x-y|)$ stands for the cost of managing the interactions between services located at $x$ and at $y$. This new choice for $G$ is more concerned with the positions of the services, and not only with the size of each pole.

These two choices and other possible models give different interesting results when one looks at the minimizers. In the first case several atoms occur in $\nu$, and $\mu$ is concentrated on balls around these centers, which corresponds to sub-cities; in the other a single-center city is obtained. The mathematical properties which are obtainable thanks to what we know from the theory of optimal transport are remarkable.

As a simple example, we will mention that the solutions of
$$(P_\nu)\quad\min_\mu W_p^p(\mu,\nu)+F(\mu);\quad \;\mbox{ for fixed } \nu\in\pical(\Omega) $$
are carachterized by
$$\mu=u\cdot\lcal^n;\quad u=(f')^{-1}\big(const-\psi_{\mu,\nu}\big)_+$$
where $\psi_{\mu,\nu}$ is a Kantorovitch potential for the transport from $\mu$ to $\nu$ and the cost $c(x,y)=|x-y|^p$.

Moreover, in the whole minimization with respect to $\mu$ and $\nu$, in the particular case  $T(\mu,\nu)= W_2^2(\mu,\nu)$ and $G(\nu)=\lambda\int_{\Omega\times\Omega}|x-y|^2\nu(dx)\nu(dy),$ $F(\mu)=||\mu||_{L^2(\Omega)}^2,$ any pair of minimizers $(\mu,\nu)$ is shaped as follows:
\begin{itemize}
\item $\mu$ is concentrated on a ball $B(x_0,r_{\lambda})$ (intersected with $\Omega$) and has a density $u$ given by $$u(x)=\frac{\lambda}{2\lambda+1}(r_{\lambda}^2-|x-x_0|^2);$$
\item $\nu$ is concentrated on the ball $B(x_0,r_{\lambda}/(2\lambda+1))$ and it is the image of $\mu$ under the homothety of ratio $(2\lambda+1)^{-1}$ and centre $x_0$;
\item $x_0$ is the barycentre of both $\mu$ and $\nu$.
\end{itemize}

The main tool for all this results is the following computation: if $\mu_\ve=(1-\ve)\mu+\ve\mu_1$, then
$$\lim_{\ve\to 0}\frac{W_p^p(\mu_\ve,\nu)-W_p^p(\mu,\nu)}{\ve}=\int \psi_{\mu,\nu}\,d(\mu_1-\mu),$$
where $\psi_{\mu,\nu}$ is, again, a Kantorovitch potential for the transport from $\mu$ to $\nu$ and the cost $|x-y|^p$. This formula says that the Kantorovitch potentials stand for Gateaux derivatives of the functional $W_p^p(\cdot,\nu)$. Thanks to standard convex analysis, it is not difficult to guess it, and to apply it to variational problems, if one thinks that the duality formula $W_p^p(\mu,\nu)=\sup \int\psi d\mu+\psi^cd\nu$ exactly says that this functional is convex and the optimal $\psi$ are the element of its subdifferential.

Notice that this kind of technique for finding optimality conditions of problem such as $(P_\nu)$ is useful in other contexts as well, and in particular, for $p=2$, when {\it gradient flows} for functional on the Wasserstein space $\W_2$ are concerned. Actually, the standard {\it minimizing movement} procedure for the gradient flow of a functional $F$ passes through discrete minimization steps for quantity like 
$$\mu\mapsto\frac{W_2^2(\mu,\nu)}{2\tau}+F(\mu),$$
where $\tau$ is a time step and a discrete sequence $(\mu_k)_k$ is built taking for $\mu_{k+1}$ the solution of $(P_{\mu_k})$. 

We finish the section by saying that other models with different costs, for instance when $T$ is no more a Wasserstein distance but comes from a congested or branched transport problem (see Section 3), have been investigated as well (see \cite{CarSan}).
 
 \section{Equilibrium structure of a city}
 
 This second part is devoted to a much more detailed model on the structure of a city which looks at an equilibrium configuration for the behavior of residents, firms and landowners. This has a much more economical taste and it has been studied by G. Carlier and I. Ekeland in \cite{carlierekeland, carlierekeland2}.
 
 The elements in this description of the city are the following:
 \begin{itemize}
 \item a domain $\Omega\subset \R^d$ which stands for the urban region we consider
 \item a measure $\mu=N(x)dx$ on $\Omega$ standing for the residents: this is unknown as well as its mass
 \item a measure $\nu=n(x)dx$ standing for jobs, which is unknown too
 \item a transportation cost $c(x,y)$ for commuting inside $\Omega$: this is given
 \item a wage function $\psi:\Omega\to\R$, where $\psi(x)$ stands for the salary that workers employed by the firm located at $x$ receive from the firm (this is an unknown of the problem)
 \item a revenue function $\phi:\Omega\to\R$ (unknown) standing for the revenues net of commuting cost that residents earn: people living at $x$ will solve $\max_y \psi(y)-c(x,y):=\phi(x)$ so as to choose where to work according to this optimization problem and, conversely, firms located at $y$ will decide whom to hire solving $\min_x \phi(x)+c(x,y)$ and getting again $\psi(y)$ as a (minimal) wage to be assured at $y$ so that there are workers who do accept to work at $y$
 \item a same utility function for all the residents $U(C,S)$ depending on their consumption level $C$ and on the quantity $S$ of land they use, as well as a fixed utility level $\bar{u}$ that every agents wants to realize: these are given as exogenous (fixed) and $\bar{u}$ may be thought as the utility level outside $\Omega$, i.e. the utility realized if one decides to move out of the city
 \item a price for residential rent $Q(x)$: at every point $x$ the residents want to choose a consumption $C$ and a land surface $S$ so that they obtain at least the utility $\bar{u}$, i.e., if $Q(x)$ is known, they solve $\min\{C+Q(x)S\,:\,U(C,S)\geq \bar{u}\}$ and they get the minimal amount of money they need. At the equilibrium this amount will necessarily be $\phi(x)$ (i.e. the money they actually have). This gives a relation between $Q$ and $\phi$ and finds the optimal value $S(x)$ as well. One obviously has $N(x)=1/S(x)$
 \item a productivity function $z:\Omega\to\R$ which is supposed to depend increasingly on $\nu$ (say, $z(x)=\nu(B(x,r))$ or $z$ is obtained through a more general convolution of $\nu$: the idea is that the productivity is higher where there is a higher concentration of workers)
 \item a production $f(z,n)$ which gives the output of a firm employing $n$ workers in a zone where the productivity is $z$
 \item a price for industrial rent $q(x)$ which is obtained, at the equilibrium, by imposing that all the surplus of the firm may be absorbed by the landlord, so that $q(x)=\max_n f(z(x),n)-\psi(x)n$. This also gives the optimal value $n(x)$.
 \end{itemize}
 
 An equilibrium is given by a pair of measures $(\mu,\nu)$, some continuous functions $z,\phi,\psi$ and a transport plan $\gamma\in\Pi(\mu,\nu)$ such that
  \begin{itemize}
 \item $\mu$ and $\nu$ have the same mass
  \item $\gamma$ and the pair $(\phi,\psi)$ are compatible in the sense that $\gamma$ is concentrated on the set $\{(x,y)\,:\,\phi(x)=\psi(y)-c(x,y)\}$ and the inequality $\phi(x)\geq\psi(y)-c(x,y)$ holds for every $(x,y)$
  \item $z$ is obtained from $\nu$ through the productivity relation that we mentioned above
 \item once $Q$ and $q$ are computed (depending on $\phi$, $\psi$ and $z$) one finds the optimal $N$ and $n$ to be equal to the densities of $\mu$ and $\nu$, respectively
\item $\mu$ is concentrated on $Q\geq q$ and $\nu$ on $q\geq Q$ (this depends on the landlords' behavior: they would not rent to residents if renting to firms is more profitable nor viceversa). 
\end{itemize}

The application of the optimal transport theory is straightforward and it allows to pose the problem as a fixed-point issue on $\mu$ and $\nu$: once $\mu$ and $\nu$ are given, one only needs to take for $\gamma$ the optimal transport plan for the cost $c$ and for  $\phi$ and $\psi$ the Kantorovitch potentials. This is what Carlier and Ekeland did, proving well-posedness results in a framework which was much more general than what was studied before in the literature (mainly one-dimensional or radially symmetric cases).
 \section{Application to Economics: Kantorovitch potential as prices or utilities}
 
 \subsection{Hotelling}
The Hotelling problem is a double-step equilibrium problem for the strategic location of $N$ firms trying to maximize their incomes from a given distribution $\mu$ of consumers in a domain $\Omega$, according to the following criterion.  Notice that the domain may be interpreted in a geographical way, or represent the different features of the goods the firms sell. 

If we know the positions $x_i$ of the firms and the prices $p_i$ that they chose, the consumer locatd at $x$ will chose where to buy his good by minimizing the sum $c(x,x_i)+p_i$ over $i=1,\dots, N$ (the cost $c(x,y)$ representing for instance the distance from $x$ to $y$ or taking into account the utility that $x$ has when he buys a product of type $y$). In this way some influence regions
$$A_i=\{x\,:\, x_i\mbox{ minimizes }c(x,x_i)+p_i\}$$
and some demands $d_i=\mu(A_i)$ are obtained. Every firm wants to maximize the profit $p_id_i$ and a {\it Nash Equilibrium} configuration for prices is a choice of the $N$ prices so that no firm wants to change its mind (i.e. changing its price $p_i$, supposing that all the other do not change their own prices). Supposing that, for every configuration of the positions of the firms, there is a unique equilibrium, every firm knows the function associating the profits to positions. An equilibrium configuration is hence a configuration where no firm wants to move in order to enhance its profit, provided the other do not move (once again, a Nash Equilibrium). The Hotelling problem exactly looks at finding such an equilibrium (see \cite{hotelling}). 

An easy but interesting link with optimal transport is the following and concerns the first step (i.e. price equilibria). The idea is: instead of taking the prices $p_i$, look at the demands $d_i$. It will be possible to reconstruct the $p_i$ from the $d_i$: in order to do that, just consider the measure $\nu=\sum_{i=1}^N d_i\delta_{x_i}$
 and prove that the function $p:\{x_1,\dots,x_N\}\to \R$ is a Kantorovitch potential for the transport from $\nu$ to $\mu$ for the cost $c$. Once this is done, the problem may be translated into a condition on $\nu$ which involves its Kantorovitch potential. Notice that writing down the precise conditions on $\nu$ involves understanding how the Kantorovitch potential depends on $\nu$, which is a very delicate issue that we will meet again.
  
\subsection{Rochet-Choné}
 
There are different models on the prices that a monopolist firm may impose for the goods it produces. One of the mathematically most interesting is the Rochet-Choné model (see \cite{RocCho}), which is an optimisation
problem under convexity constraint. The convex structure comes from the simplifying assumption that the space of goods $y$ and the space of consumers $x$ are subsets of $\R^N$ and they are coupled through the function $(x,y)\mapsto x\cdot y$ representing the utility that a consumer of type $x$ has in buying $y$. Once the distribution $\mu$ of consumers is known, the firm may choose the price for its good, i.e. a function $p:Y\to [0,\infty[$, defined on the goods space $Y$; then every consumer $x$ choses what to buy by solving
$$\max_y\; x\cdot y -p(y)$$
and getting a utility $u(x):=\max_y x\cdot y -p(y)$, realized by a good $y_x$. The firm may reconstruct its total gain by integrating $p(y(x))-C(y(x))$ (if $C(y)$ is the cost for producing $y$) . The total profit is hence given by
$$\int_X \left(p(y(x))-C(y(x))\right)\mu(dx)=\int_X \left(x\cdot y(x)-u(x)-C(y(x))\right)\mu(dx).$$
One can also notice that $y_x=\nabla u(x)$ (differentiating the expression of $u$) and hence the maximization of the profit is a problem that may be stated in terms of $u$
$$\max F(u)=\int_X \left(x\cdot \nabla u(x)-u(x)-C(\nabla u(x))\right)\mu(dx)$$
where the constraint on $u$ are convexity (from its defintion) and positivity ($u\geq 0$ is a consequence of the fact that consumers do not buy if they get a negative utility: it may be stated saying that a certain ``empty'' good called $0$ belongs to $Y$ and that we impose $p(0)=0$; the firm is not allowed to charge for buying this empty good but this good interests nobody) and a constraint on the gradient: $\nabla u\in Y$. This is the minimization problem under convexity constraint we referred to. It falls into the framework of the convexity-constrained problems studied for instance by Carlier and Lachand-Robert  (see \cite{CarLac}), where some $C^1$ regularity results are also proven. The same class of problems also includes the well-known Newton Problem of minimal resistance. For both the problems, some numerical insights in particular cases exist, but lots of information lack.

An interesting change of variable, using the image measure $\nu=(\nabla u)_\sharp \rho$, is possible, since every measure is the image of $\rho$ through the gradient of a convex function (which is exactly the well known result by Brenier in transport theory, see \cite{polarization}). It is interesting to link this reformulation to optimal transport.The most natural cost to be considered would be the scalar product but we know that considering $-x\cdot y$ or $\frac 12 |x-y|^2$ is the same. Hence, we may rewrite the previous problem as
$$\min_{u\;\mbox{ convex}} \tilde{F}(u)=\int_X \left(\frac{|x|^2}{2}-x\cdot \nabla u(x)+\frac{|\nabla u|^2}{2}+u(x)+\tilde{C}(\nabla u(x)\right)\mu(dx)$$
where $\tilde{C}(z)=C(z)-|z|^2/2$ and we are allowed to add the term in $|x|^2/2$ since it does not depend on $u$. 

We can in the end rewrite the problem in terms of $\nu$ as
$$\min_{\nu\in\pical(Y)} G(\nu)=\frac 12 W_2^2(\mu,\nu)+ \int_Y\tilde{C}d\nu + \int_X u_\nu \,d\mu,$$
$u_\nu$ being for a measure $\nu$ the unique convex function satisfying $\nabla u_\#\mu=\nu$ and $\min u= 0$ (which is obtained as a Kantorovitch potential for the cost $-x\cdot y$ or $\frac 12 |x-y|^2$).

This kind of functional may be considered via the transport theory. Existence of a minimizer is easy and the interesting point is finding optimality conditions. The difficult part is handling the term
$$\nu\mapsto\int_X u_\nu \,d\mu.$$
For getting optimality conditions, it would be useful to differentiate this term with respect to variations of $\nu$. Yet, computing
$$\lim_{\ve\to 0}\frac{u_{\nu_\ve}-u_\nu}{\ve},\qquad \nu_\ve=\nu+\ve(\tilde{\nu}-\nu)$$
is a challenging issue; 
 possible strategies include the linearisation of the Monge-Ampère equation but lots of questions are open.

\end{document}